\documentclass[11pt]{article}
\usepackage{xcolor}
\usepackage{amsmath}
\usepackage{amsfonts}
\usepackage{amsthm}
\usepackage{enumerate}
\usepackage[english]{babel}
\usepackage{calrsfs}
\newcommand{\KK}{\mathbb K}

\newcommand{\cP}{\mathcal P}

\newcommand{\cS}{\mathcal S}

\newcommand{\cL}{\mathcal L}
\newcommand{\cG}{\mathcal G}

\newcommand{\cLh}{{\mathcal L}_h}

\newcommand{\PG}{\mathrm{PG}}

\newcommand{\fL}{\mathfrak L}

\newcommand{\rk}{\mathrm{rank}}

\newcommand{\ve}{\varepsilon}

\newtheorem{theorem}{Theorem}[section]
\newtheorem{lemma}[theorem]{Lemma}
\newtheorem{corollary}[theorem]{Corollary}
\newtheorem{prop}[theorem]{Proposition}
\newtheorem{res}[theorem]{Result}
\newtheorem{example}[theorem]{Example}
\theoremstyle{definition}
\newtheorem{remark}{Remark}

\makeatletter
\let\ltxxlabel\ltx@label
\makeatother

\title{Characterizations of symplectic polar spaces}
\author{I. Cardinali, H. Cuypers, L. Giuzzi and A. Pasini}
\begin{document}

\maketitle

\begin{abstract}
\noindent
A polar space $\cS$ is called symplectic if it admits a projective embedding $\ve\colon \cS \rightarrow \PG(V)$ such that the image $\ve(\cS)$ of $\cS$ by $\ve$ is defined by an alternating form of $V$. In this paper we characterize symplectic polar spaces in terms of their incidence properties, with no mention of peculiar properties of their embeddings. This is relevant especially when $\cS$ admits different (non isomorphic) embeddings, as it is the case when $\cS$ is defined over a field of characteristic $2$.
\end{abstract}

\leftline{{\bfseries MSC:} 51A50, 51B25, 51E24 }
\leftline{{\bfseries Keywords:} polar spaces, embeddings, hyperplanes, hyperbolic lines }

\section{Introduction}

We assume the reader is familiar with basics on polar spaces, for which we refer to \cite{BC}, \cite{S} and \cite{T}. So, we are not going to recall the definition of polar space, singular subspaces, (geometric) hyperplanes and, more generally, subspaces of a polar space. Basic definitions of embeddings of polar spaces can be found in Section \ref{emb}.
By assumption, all polar spaces to be considered in this paper are non-degenerate, thick-lined (all of their lines have at least three points) and do admit lines, namely they have rank at least 2. We do not assume the rank to be finite; we shall need this assumption only in a few results. 

We adopt the following notation. The symbol $\perp$ stands for collinearity, with the convention that every point is collinear with itself. Given a polar space $\cS$, let $P$ be its set of points. For $p\in P$ we put $p^\perp := \{x\in P~|~x\perp p\}$ and, for $X\subseteq P$, we put $X^\perp := \cap_{x\in X}x^\perp$, also writing $\{X, Y\}^\perp$ for $(X\cup Y)^\perp$, $X^{\perp\perp}$ for $(X^\perp)^\perp$ and so on. With $X$ as above, we denote by $\langle X\rangle$ the subspace of $\cS$ spanned by $X$.

\subsection{Two known characterizations of symplectic polar spaces}\label{prel}
 
In a strict sense, a {\em symplectic polar space} is the full subgeometry of the projective geometry $\PG(V)$ of a vector space $V$, with the same points as $\PG(V)$ and the lines of which are the lines of $\PG(V)$ which are totally isotropic for a given non-degenerate alternating form of $V$. In a less strict sense, a polar space $\cS$ is called {\em symplectic} if it admits an embedding $\ve:\cS\rightarrow\PG(V)$ such that the embedded polar space $\ve(\cS) \subseteq \PG(V)$ is symplectic in the previous sense; equivalently, $\ve$ maps the point-set of $\cS$ surjectively onto the point-set of $\PG(V)$. When $\cS$ is defined over a field of characteristic different from 2 these two definitions essentially coincide but in characteristic 2 they are different. Indeed,  every  polar space defined over a division ring of characteristic different from $2$ admits a unique embedding (unique modulo isomorphisms, of course), except for two particular families of generalized quadrangles, which however have nothing to do with symplectic spaces. See for example \cite{BC,S,J2,T}. In contrast, in characteristic $2$, a symplectic polar space always admits different embeddings. In the sequel we follow both definitions: when saying that an embedded polar space $\ve(\cS)$ is symplectic we cling to the first definition while, when regarding $\cS$ as an abstract object, we hold to the latter.      

We say that a hyperplane $H$ of $\cS$ is {\em singular} with {\em deepest point $p$} if $H = p^\perp$ for a (necessarily unique) point $p$ of $\cS$. When $\cS$ is embeddable, all singular hyperplanes of $\cS$ arise from each of the embeddings of $\cS$. The following is folklore.
 
\begin{prop}\label{res 1}
Let $\ve$ be an embedding of $\cS$. Then $\ve(\cS)$ is symplectic of finite rank if and only if the hyperplanes of $\cS$ which arise from $\ve$ are precisely the singular ones.
\end{prop}

As noticed above, every embeddable polar space of rank at least 3 defined over a division ring of characteristic different
from $2$ admits a unique embedding and, as follows by Shult \cite[7.4.8]{S} or \cite[Proposition 5.2]{CS90}, all its hyperplanes arise from that embedding. 
Accordingly, Proposition \ref{res 1} implies the following.

\begin{prop}\label{res 1 bis}
Let $\cS$ be embeddable and of rank at least $3$. Let $\KK$ be the underlying division ring of $\cS$. All hyperplanes of $\cS$ are singular if and only if $\cS$ is symplectic, of finite rank and $\mathrm{char}(\KK) \neq 2$.
\end{prop}

\begin{remark}
The statement  $\rk(\cS) < \infty$ cannot be dropped from Proposition \ref{res 1} and Proposition \ref{res 1 bis}. Indeed let $\ve(\cS)\subseteq \PG(V)$ be symplectic. Proposition \ref{res 1} claims for the existence of a linear bijection between $V$ and its dual $V^*$, but no such bijection exists when $\rk(\cS)$ is infinite. Indeed if this is the case then $\dim(V)$ is also infinite and therefore $\dim(V) < \dim(V^*)$.
\end{remark} 

\begin{remark}
The hypothesis that $\cS$ is embeddable cannot be dropped from Proposition \ref{res 1 bis}. Indeed, as proved by Cohen and Shult \cite{CS90}, all hyperplanes of a non-embeddable polar space of rank 3 are singular. 
\end{remark}

\begin{remark}\label{cam}
In Proposition \ref{res 1 bis} we have assumed that $\rk(\cS) > 2$. Suppose that $\rk(\cS) = 2$. If $\cS$ is finite with lines of size $q+1$, for an odd integer $q > 1$, then $\cS$ is symplectic if and only if all of its hyperplanes are singular (see Payne and Thas~\cite{PT}), like in the statement of Proposition \ref{res 1 bis}. On the other hand, let $\cS$ be symplectic and infinite. Then $\cS$ admits infinitely many ovoids (Cameron \cite{Cam}), most of which are non-classical, namely do not arise from any embedding of $\cS$. (In fact when $\cS$ is defined over a field of characteristic other than 2 none of them is classical.) So, the statement of Proposition \ref{res 1 bis} fails to hold when $\rk(\cS) = 2$ and $\cS$ is infinite. 
\end{remark}

In view of the next characterizations of symplectic polar spaces, we recall the definition of hyperbolic lines. Let $a$ and $b$ be two non-collinear points of $\cS$. Then $\{a,b\}^{\perp\perp}$ is a set of mutually non-collinear points called a \emph{hyperbolic line} through $a$ and $b$. 

Let $P$ and $\cL$ be respectively the point-set and the set of lines of $\cS$; denote by $\cLh$ the set of all hyperbolic lines of $\cS$ and put $\fL(\cS) := (\cP,\cL\cup\cLh)$. It is well known that $\fL(\cS)$ is a linear space (any two points are joined by a unique line). Obviously, $\cS$ is a full subgeometry of $\fL(\cS)$, with exactly the same points as $\fL(\cS)$. The next statement is folklore.

\begin{prop}\label{main thm 1}
The polar space $\cS$ is symplectic if and only if $\fL(\cS)$ is a projective space.
\end{prop}

\subsection{Main results}\label{Main sec} 
Our aim in this paper is to obtain `embedding-free' characterizations of symplectic polar spaces, namely characterizations where no property
of any embedding is mentioned but possibly the assumption that the space is embeddable. Proposition \ref{main thm 1} seemingly meets these requirements but actually it does not. Indeed if $\fL(\cS)$ is a projective space then the identity mapping $p\in \cS \mapsto p\in \fL(\cS)$ is just the symplectic embedding of $\cS$. Thus Proposition \ref{main thm 1} amounts to the following trivial embedding-committed characterization: $\cS$ is symplectic if and only if it admits an embedding $\ve:\cS\rightarrow\PG(V)$ where $\ve$ maps the point-set of $\cS$ onto the point-set of $\PG(V)$. Proposition \ref{res 1 bis} looks better in the perspective we have chosen for this paper, but it misses the characteristic 2 case and fails to hold when the rank is either 2 or infinite. The next theorem, to be proved in Section \ref{MT1}, improves Proposition \ref{res 1 bis}. It holds with no restrictions on either the rank of $\cS$ or the characteristic of its underlying division ring.   

\begin{theorem}\label{mainthm 3}
An embeddable polar space $\cS$ is symplectic if and only if the following holds:
\begin{itemize}
\item[$(\ast)$] if a hyperplane $H$ of $\cS$ contains the trace $\{a,b\}^\perp$ of two non-collinear points $a$ and $b$ of $\cS$ then $H$ is singular.
\end{itemize}
\end{theorem}
Maximal singular subspaces are often called {\em generators} in the literature. We shall also adopt that terminology here, with the warning that in general a polar space of infinite rank admits generators of different dimensions (see e.g. \cite{Pas}). The following Lemma will be proved in Section \ref{MT2}.  

\begin{lemma}\label{nuovo}
Let $\cS$ be embeddable. Then $\cS$ satisfies property $(\ast)$ of Theorem {\rm \ref{mainthm 3}} if and only if both the following hold in $\cS$ for any two non-collinear points $a$ and $b$:
\begin{itemize}
\item[$(A)$] if $M$ and $N$ are respectively a generator of $\cS$ and a generator of the polar space $\{a,b\}^\perp$ and $M \supseteq N$, then $M\cap\{a,b\}^{\perp\perp}\neq\emptyset$;
\item[$(B)$] if a hyperplane $H$ of $\cS$ contains $\{a,b\}^{\perp}$, then $H$ contains a generator $M$ of $\cS$ such that $M\cap\{a,b\}^\perp$ is a generator of $\{a,b\}^\perp$. 
\end{itemize}
\end{lemma}
Consequently:

\begin{corollary}\label{main thm 2}
An embeddable polar space is symplectic if and only if both properties $(A)$ and $(B)$ of Lemma {\rm \ref{nuovo}} hold in it for any pair of non-collinear points. 
\end{corollary}

\begin{remark}
The hypothesis that $\cS$ is embeddable cannot be dropped from any of Theorem \ref{mainthm 3}, Lemma \ref{nuovo} or Corollary \ref{main thm 2}. Indeed non-embeddable polar spaces exist which satisfy $(\ast)$; hence they also satisfy $(B)$, since $(\ast)$ implies $(B)$ regardless of embeddability, as one can see by recalling that non-embeddable polar spaces have rank at most 3 and rephrasing $(B)$ in the finite rank case. For instance, $(\ast)$ trivially holds in every non-embeddable polar space of rank 3, since all hyperplanes of such a polar space are singular (Cohen and Shult \cite{CS90}). Non-embeddable polar spaces of rank 3 also satisfy property $(A)$ of Lemma \ref{nuovo} (Section \ref{MT3}, Result \ref{MTA}). Non-embeddable generalized quadrangles also exist which satisfy $(\ast)$ but not $(A)$. Two examples of this kind will be discussed in Section \ref{MT3}. 
\end{remark} 

As proved in \cite[Lemma 5.5]{PVM} (see also this paper, Proposition \ref{A0}), when $\rk(\cS) < \infty$ a pair $\{a,b\}$ of non-collinear points satisfies $(A)$ if and only if $\{a,b\}^{\perp\perp} = \{N, N'\}^\perp$ for any choice of mutually disjoint generators $N$ and $N'$ of $\{a,b\}^\perp$. We say that a pair $\{a,b\}$ with this property is {\em regular}. Indeed the above property is a natural generalization of regularity as defined for generalized quadrangles (Van Maldeghem \cite{HVM1}, also Payne and Thas \cite{PT} for the finite case). On the other hand, without assuming that $\rk(\cS) < \infty$, let $\cS$ be embeddable. Then, as we shall prove in Section \ref{B centric}, property $(B)$ of Lemma \ref{nuovo} is equivalent to the following: for any point $c$ collinear with neither $a$ nor $b$ the triad $\{a,b,c\}$ is {\em centric}, namely $\{a,b,c\}^\perp$ contains a generator of $\{a,b\}^\perp$. (This is also a quite natural generalization of the concept of being centric as defined for generalized quadrangles.) Accordingly, in the finite rank case Corollary \ref{main thm 2} can be rephrased in the following way, which reminds us of well known characterizations of symplectic generalized quadrangles. 

\begin{corollary}\label{main thm 2bis}
An embeddable polar space $\cS$ of finite rank is symplectic if and only if all pairs of non-collinear points of $\cS$ are regular and all triads of mutually non-collinear points of $\cS$ are centric. 
\end{corollary}

Referring the reader to Section \ref{sect5} for more on properties $(A)$ and $(B)$ we now turn to our second main theorem. Its statement has been suggested by Hendrik Van Maldeghem \cite{HVM2} to the fourth author of this paper. 

\begin{theorem}\label{Hendrik 1}
A polar space $\cS$ is symplectic if and only if every singular hyperplane of $\cS$ meets every hyperbolic line of $\cS$ non-trivially.
\end{theorem}
Note that in this theorem we do not assume $\cS$ to be embeddable. A proof of Theorem \ref{Hendrik 1} will be given in Section \ref{sect6}.

Obviously, all subspaces of the linear space $\mathfrak{L}(\cS) = (P, \cL\cup\cL_h)$ introduced at the end of Section~\ref{prel} are subspaces of $\cS$; in particular, all geometric hyperplanes of $\mathfrak{L}(\cS)$ are hyperplanes of $\cS$. Conversely, the singular hyperplanes of $\cS$ are closed under taking hyperbolic lines, hence they are subspaces (in fact maximal subspaces) of $\mathfrak{L}(\cS)$. However in general they are not hyperplanes of $\mathfrak{L}(\cS)$. Indeed, according to Theorem \ref{Hendrik 1}, all singular hyperplanes are hyperplanes of $\mathfrak{L}(\cS)$ if and only if $\cS$ is symplectic. Turning to non-singular hyperplanes, the following will be proved in Section \ref{sect6}:

\begin{prop}\label{nuovo H}
Let $\cS$ be a symplectic polar space of rank at least $2$. A hyperplane $H$ of $\cS$ arises from the symplectic embedding of $\cS$ if and only if it contains a hyperbolic line.
\end{prop}

Every hyperplane of a symplectic polar space $\cS$ contains a pair of non-collinear points. Consequently, with $\cS$ and $H$ as in Proposition \ref{nuovo H}, hyperbolic lines exist which meet $H$ in at least two points but $H$ contains none of them; hence $H$ is not even a subspace of $\mathfrak{L}(\cS)$. Nothing like this occurs when the underlying field of $\cS$ has characteristic different from $2$ and either $2 <\rk(\cS) < \infty$ or $\rk(\cS) = 2$ but $\cS$ is finite, since in each of these two cases all hyperplanes of $\cS$ are singular (Proposition \ref{res 1 bis} and Remark \ref{cam}). So, Theorem \ref{Hendrik 1} implies the following:

\begin{corollary}\label{Hendrik 2}
Let $\cS$ be a polar space.
Then all hyperplanes of $\cS$ are hyperplanes of $\mathfrak{L}(\cS)$ if and only if $\cS$ is 
symplectic and defined over a field of characteristic different from $2$, with finite rank at least $3$, or,  with rank $2$ and being finite. 
\end{corollary}

\begin{remark}\label{H rem} 
When $\cS$ is symplectic of rank $2$, hyperbolic lines and traces of pairs of non-collinear points are the same; explicitly, if $a$ and $b$ are non-collinear points of $\cS$ and $c$ and $d$ distinct (necessarily non-collinear) points of $\{a,b\}^\perp$, then $\{c,d\}^\perp = \{a,b\}^{\perp\perp}$ (in short, $\cS$ is regular). Therefore, by property $(\ast)$ of Theorem \ref{mainthm 3}, no ovoid of $\cS$ contains a hyperbolic line. However, as noticed in Remark \ref{cam}, if moreover $\cS$ is infinite then $\cS$ admits ovoids. So, the hypothesis that $\cS$ is finite when $\rk(\cS) = 2$ cannot be dropped from Corollary~\ref{Hendrik 2}.
\end{remark}

\noindent
\textbf{Organization of the paper.} In Section \ref{emb} we recall all we need from the theory of embeddings of polar spaces. Sections \ref{Proof MT} and \ref{sect6} contain the proofs of Theorems \ref{mainthm 3} and Lemma \ref{nuovo} and, respectively, Theorem \ref{Hendrik 1} and Proposition \ref{nuovo H}. In Section \ref{sect5} we discuss properties $(A)$ and $(B)$ of Lemma \ref{nuovo}, focusing on properties equivalent to them.  

\section{Embeddings of polar spaces}\label{emb}

Suppose $\cS=(P,\mathcal{L})$ is a polar space and $\PG(V)$ a projective geometry of a vector space $V$ defined over some skew field.
Then an \emph{embedding} $\ve:\cS\rightarrow \PG(V)$ is an injective map from $P$ to the point set of $\PG(V)$ with $\ve(P)$ generating $\PG(V)$ and mapping each line $\ell\in \mathcal{L}$ surjectively to a line of $\PG(V)$.
(Here we consider lines of both $\cS$ and $\PG(V)$ to be sets of points.)

Suppose $\ve:\cS\rightarrow \PG(V)$ is an embedding of the polar space $\cS$.
For any subset $X$ of points of $\PG(V)$ we denote by $\langle X\rangle_V$ the span of $X$ in $\PG(V)$.
Moreover, by  $\ve^{-1}(X)$ we denote the set of all points in $P$ mapped by $\ve$ to points in $X$, i.e.  $\ve^{-1}(X\cap \ve(P))$.

Embeddings of polar spaces have been studied extensively. See for example \cite{BC} or \cite{S} and the reference therein.

The properties of embeddings important to us are given below.

\begin{prop}\label{hyperplane embedding}
Let $x\in P$. Then  $\langle \ve(x^\perp)\rangle_V$ is a hyperplane of $\PG(V)$
with  $\ve^{-1}(\langle \ve(x^\perp)\rangle_V)=x^\perp$.
\end{prop}
\begin{proof}
See for example \cite[7.7.4]{S} or \cite[5.4]{J2}.\hfill  \end{proof}

If $x\in P$, then by $\ve(x)^{{\perp_\ve}}$ we denote the hyperplane $\langle \ve(x^\perp)\rangle_V$.
For any subset $X$ of $\PG(V)$ we denote by $X^{\perp_{\ve}}$ the
intersection of all hyperplanes $\ve(x)^{{\perp_\ve}}$ with $x\in \ve^{-1}(\langle X\rangle_V)$.

\begin{lemma}\label{intersection of hyperplanes}
If $x,y\in P$ are two  points, then
$$\{x,y\}^\perp=\ve^{-1}(\ve(x)^{{\perp_\ve}}\cap \ve(y)^{{\perp_\ve}})=\ve^{-1}(\{\ve(x),\ve(y)\}^{{\perp_\ve}}).$$
\end{lemma}
\begin{proof}
Clearly a point $z\in \{x,y\}^\perp$ maps to a point of $\ve(x)^{{\perp_\ve}}\cap \ve(y)^{{\perp_\ve}}$. Moreover, if $z\in P$ maps to a point in $\ve(x)^{{\perp_\ve}}\cap \ve(y)^{{\perp_\ve}}$, then
by Proposition \ref{hyperplane embedding}, it is in $x^\perp$ as well as in $y^\perp$, and thus in $\{x,y\}^\perp$. This proves the first equality.

Now $\{\ve(x),\ve(y)\}^{{\perp_\ve}}$ is the intersection of all hyperplanes $\ve(z)^{{\perp_\ve}}$, where $z$ is a point which maps
to a point $\ve(z)$ on the projective line on $\ve(x)$ and $\ve(y)$. Clearly $\{\ve(x),\ve(y)\}^{{\perp_\ve}}\subseteq \ve(x)^{{\perp_\ve}}\cap \ve(y)^{{\perp_\ve}}$. If $u\in P$ maps to a point in $\ve(x)^{{\perp_\ve}}\cap \ve(y)^{{\perp_\ve}}$,
then it is collinear to both $x,y$. So, $\ve(u)^{{\perp_\ve}}$
contains $\ve(x)$ and $\ve(y)$ and then also $\ve(z)$. But then $z\perp u$ and $\ve(u)\in \ve(z)^{{\perp_\ve}}$.
This proves  $\{\ve(x),\ve(y)\}^{{\perp_\ve}}\supseteq \ve(x)^{{\perp_\ve}}\cap \ve(y)^{{\perp_\ve}}$ and hence the second equality.\hfill \end{proof}

\noindent
\textbf{Hyperplanes and embeddings.}
Let $\cS$ be an embeddable polar space and let $\ve:\cS\rightarrow\PG(V)$ be an embedding of $\cS$. A subspace $X$ of $\cS$ {\em arises from} $\ve$ if $\ve^{-1}(\langle \ve(X)\rangle_V) = X$. In particular, a geometric hyperplane $H$ of $\cS$ arises from $\ve$ if and only if $\ve(H)$ spans a projective hyperplane of $\PG(V)$. 

Proposition \ref{hyperplane embedding} implies that all singular hyperplanes of $\cS$ arise from any embedding of $\cS$. 

\bigskip

\noindent
{\textbf{Embeddability.} All polar spaces of rank at least 4 admit an embedding (Tits \cite[Chapter 8]{T}, Buekenhout and Cohen \cite[Chapter 8]{BC}, also Cuypers, Johnson and Pasini \cite{CJP}); a polar space $\cS$ of rank 3 is embeddable if and only if its planes are desarguesian but $\cS$ is not the line-grassmannian of a 3-dimensional projective geometry over a non-commutative division ring (Tits \cite[Chapter 8]{T}, Buekenhout and Cohen \cite[Chapter 8]{BC}). Exactly one family of (necessarily non-embeddable) polar spaces of rank 3 with non-desarguesian planes exist, independently discovered by Freudenthal and Tits in the fifties of last century and described in \cite[Chapter 9]{T} (see also M\"{u}hlherr \cite{M} and De Bruyn and Van Maldeghem \cite{DBVM1, DBVM2}). For further reference, we shall call them {\em Freudenthal-Tits} polar spaces.} \\

For almost all embeddable polar  spaces defined over a division ring of characteristic not $2$, there is
essentially one unique embedding. Indeed, according to Tits \cite[chp. 8]{T} (see also Johnson \cite{J1, J2}), this holds for allowed polar spaces $\cS$, except precisely when $\cS$ is a grid with lines of size at least 5 or a generalized quadrangle as described in \cite[8.6 {$\mathrm{(II)(a)}$}]{T}; the latter exceptional generalized quadrangles are defined over quaternion division rings and admit two non-isomorphic embeddings, both of which have dimension 4. 

In case the polar space is defined over a division ring in characteristic $2$, there might be more embeddings.
As we are mainly concerned with symplectic spaces we provide an example based on \cite{DB}.

Let $V$ be a vector space of dimension at least $4$ over a field $\mathbb{K}$ of characteristic $2$ equipped with a non-degenerate symplectic form $f$
and basis $B$.  Let $\cS$ be the corresponding symplectic polar space.
We will construct another embedding of $\cS$.

For this, let $R$ be a $\mathbb{K}$-vector space  equipped with an anisotropic quadratic form $Q$ whose associated bilinear form is $0$.
Notice that for $r,r'\in R$ we have $Q(r)=Q(r')$ if and only if $r=r'$. Then consider $\hat{V}=V\oplus R$
with symplectic form  $\hat{f}$
defined by
$$\hat{f}(v,w)=
\begin{cases}f(v,w) & \mathrm{if} \ v,w\in V\\
\hspace{5 mm} 0 & if \ v\in R \ \mathrm{or}\ w\in R\\
\end{cases}$$
and quadratic form $\hat{Q}$ given by
$$\hat{Q}(v)=\begin{cases}1\ \mathrm{if} \ v\in B\\
                          Q(v) \ \mathrm{if} \ v\in R
                          \end{cases}$$
and for all $v,w\in \hat{V}$
$$\hat{Q}(v+w)=\hat{Q}(v)+\hat{Q}(w)+\hat{f}(v,w).$$
Let $X$ be a subspace  of $\hat{V}$ containing $R$ as a hyperplane. Fix a vector $v\in X\setminus R$. Every $1$-space of $X$ not in $R$ contains then a vector  $x=v+r$ with for some $r\in R$.
Hence  
$$\hat{Q}(x)=\hat{Q}(v)+\hat{Q}(r)= \hat{Q}(v)+Q(r).$$
In particular, there is at most one $1$-space in $X$ on which $\hat{Q}$ vanishes.
The quotient map from $\hat{V}\rightarrow \hat{V}/R\simeq V$ induces an isomorphism of the polar space defined by
$\hat{Q}$ and $\cS$ when $Q$ takes all values in $\mathbb{K}$, while its inverse map provides us with an embedding
of $\cS$ different from the symplectic embedding.  

Such form $Q$ on a vector space $R$ always exists. Indeed, take a basis $\Lambda$ of $\mathbb{K}$ over $\mathbb{K}^2$.
Then take $R$ to be the $\mathbb{K}$-space with formal basis $(b_\lambda)_{\lambda\in \Lambda}$.
Let $r\in R$, then $r$ can uniquely be expressed as $\mu_1b_{\lambda_1}+\cdots+\mu_kb_{\lambda_k}$
and we can define $Q(r)=Q(\mu_1b_{\lambda_1}+\cdots+\mu_kb_{\lambda_k})=\mu_1^2\lambda_1+\cdots +\mu_k^2\lambda_k$.
As $\Lambda$ is a basis for $\mathbb{K}$ over $\mathbb{K}^2$, we do find that $Q$ takes all values of $\mathbb{K}$.

\section{Proof of Theorem \ref{mainthm 3} and Lemma \ref{nuovo}}\label{Proof MT} 

\subsection{Proof of Theorem \ref{mainthm 3}}\label{MT1}

\textbf{The `if' part.} Assume $(\ast)$ on $\cS$ and let $\ve:\cS\rightarrow \PG(V)$ be an embedding of $\cS$. Let $V^*$ be the dual of $V$ and consider the mapping $\ve^*:\cS\rightarrow\PG(V^*)$ defined as follows: $\ve^*(p) = \ve(p)^{{\perp_\ve}}$, for every point $p\in \cS$. 

Let $a$ and $b$ be distinct points of $\cS$ and let $\overline{H}$ be a projective hyperplane of $\PG(V)$ containing $\{\ve(a),\ve(b)\}^{{\perp_\ve}}$. As $\{\ve(a),\ve(b)\}^{{\perp_\ve}}$ has codimension $2$ in $\PG(V)$, we have $\overline{H} = \langle \{\ve(a), \ve(b)\}^{{\perp_\ve}}, \bar{p}\rangle_V$ for any point $\bar{p}\in \overline{H}\setminus\{\ve(a), \ve(b)\}^{{\perp_\ve}}$. 

Suppose first that $a$ and $b$ are collinear. Then $\{\ve(a), \ve(b)\}^{{\perp_\ve}}$ contains the projective line $\bar{\ell} = \ve(\langle a, b\rangle)$ and, with $\bar{p}$ as above, $\bar{p}^{{\perp_\ve}}$ meets $\bar{\ell}$ in a point, say $\bar{c} = \ve(c)$. Therefore $\overline{H} \subseteq \ve(c)^{{\perp_\ve}}$. Hence $\overline{H} = \ve(c)^{{\perp_\ve}}$, since $\ve(c)^{{\perp_\ve}}\subset \PG(V)$.  

Suppose now that $a$ and $b$ are non-collinear. The $\ve$-preimage $H := \ve^{-1}(\overline{H})$ of $\overline{H}$ is a geometric hyperplane of $\cS$ and contains $\{a,b\}^\perp$. By $(\ast)$, we have $H = c^\perp$ for a suitable point $c$ of $\cS$. Accordingly, $\overline{H} = \ve(c)^{{\perp_\ve}}$. Clearly, $c^\perp\supseteq \{a,b\}^{\perp}$, namely $c \in \{a, b\}^{\perp\perp}$. 

So far we have proved that every hyperplane of $\PG(V)$ containing $\ve(\{a,b\}^\perp)$ belongs to the image of $\ve^*$. Moreover, if $a\perp b$ (respectively $a\not\perp b$) then $\ve^*$ maps the line $\langle a, b\rangle$ of $\cS$ (the hyperbolic line $\{a,b\}^{\perp\perp}$ of $\cS$) bijectively onto the line of $\PG(V^*)$ through $\ve(a)^{{\perp_\ve}}$ and $\ve(b)^{{\perp_\ve}}$. In the end, $\ve^*$ defines an isomorphism from the linear space $\mathfrak{L}(\cS)$ defined at the end of Section \ref{prel} to a subspace of $\PG(V^*)$. Consequently, $\mathfrak{L}(\cS)$ is a projective space. Hence $\cS$ is symplectic, by Proposition \ref{main thm 1}. \hfill $\Box$

\medskip

\noindent
\textbf{The `only if' part.} 
Suppose $\cS$ is symplectic.
Let $a$ and $b$ be non-collinear points of $\cS$ and $H$ a hyperplane of $\cS$ containing $\{a,b\}^\perp$.
Let $\ve:\cS\rightarrow\PG(V)$ be the symplectic embedding of $\cS$ with respect to some non-degenerate symplectic form $f$.
Denote the corresponding symplectic polarity by $\perp_f$.
Then $V=L{\perp_f}M$, where $M:= \langle\ve(\{a,b\}^\perp)\rangle_V$ is of codimension $2$ and  $L=\langle\ve(\{a,b\}^{\perp\perp})\rangle_V$ of dimension $2$.
As  $\{a,b\}^\perp$ is not a hyperplane of $\cS$, there is a point $c\in H$ which is not in $\{a,b\}^\perp$.
But then $\langle \ve(c),L\rangle_V$ meets $M$ in a $1$-space which equals $\ve(d)$ for some point $d$  of $\cS$.
As $\ve(d){\perp_f} \langle\ve(d), L\rangle_V$, we find $\ve(d){\perp_f} \ve(c)$ and hence $d\perp c$.
Let $e$ be the point on $\cS$ with $\ve(e)$ in the intersection of $\langle \ve(c),\ve(d)\rangle_V$ and $L$. 
Then $e^\perp$, which is generated by $c$ and $\{a,b\}^\perp$, is contained in $H$.
However, as $e^\perp$ is a maximal subspace of $\cS$,  that implies that $e^\perp=H$.   \hfill $\Box$

\subsection{Sub-generators}\label{sub-gen}

Before to turn to the proof of Lemma \ref{nuovo}, we prove a few elementary properties of generators and their hyperplanes, to be freely used in the sequel. In what follows $\cS$ is a polar space of arbitrary (possibly infinite) rank $\rk(\cS) \geq 2$. 

\begin{lemma}\label{sub-gen 1}
If a singular subspace of $\cS$ is a hyperplane of at least one of the generators of $\cS$ which contain it, then it is a hyperplane in every generator which contains it.
\end{lemma} 
\begin{proof}
Let $X$ be a singular subspace of $\cS$ and suppose that $M := \langle X, a\rangle$ is a generator for some $a\in X^\perp\setminus X$.  Let $N\neq M$ be another generator of $\cS$ containing $X$. Clearly $a \not \in N$ otherwise $M \subseteq N$, hence $M = N$ because $M$ is a generator; a contradiction with the choice of $N$. So, $M\cap N = X$. Consequently, $a^\perp\cap N$ is a hyperplane of $N$ and contains $X$. We have $a^\perp\cap N = X$. Indeed, if otherwise, let $b \in (a^\perp\cap N)\setminus X$. Then $b \perp M$, hence $b\in M$ since $M$ is a generator. However $M\cap N = X$ while $b \in N\setminus X$; contradiction.   \end{proof} 

Henceforth, the singular subspaces which occur as hyperplanes in generators of $\cS$ will be called {\em sub-generators}.    

\begin{prop}\label{sub-gen 2}
Let $a$ and $b$ be two non-collinear points of $\cS$. Then the generators of the polar space $\{a,b\}^\perp$ are precisely the sub-generators of $\cS$ contained in $\{a,b\}^\perp$. 
\end{prop}
\begin{proof} In view of Lemma \ref{sub-gen 1}, it is sufficient to prove that if $X$ is a generator of $\{a,b\}^\perp$ then $\langle X, a\rangle$ is a generator of $\cS$. Let $M$ be a generator of $\cS$ containing $\langle X, a\rangle$. Then $Y := b^\perp\cap M$ is a hyperplane of $M$ (since $b\not\perp a\in M$). Moreover, $X \subseteq Y\subseteq \{a,b\}^\perp$. However $X$ is a generator of $\{a,b\}^\perp$. Hence $Y = X$, namely $X$ is a hyperplane of $M$. Therefore $\langle X, a\rangle = M$, since $X \subset \langle X, a\rangle \subseteq M$. \end{proof} 

\subsection{Proof of Lemma \ref{nuovo}}\label{MT2}

We are now ready to prove Lemma \ref{nuovo}. We first prove the `only if' part of Lemma \ref{nuovo}, next we shall turn to the `if' part.  \\

\noindent
\textbf{The `only if' part.} Assume $(\ast)$ on $\cS$ and let $a$ and $b$ be two non-collinear points of $\cS$. Let $H$ be a hyperplane of $\cS$ containing the trace $\{a,b\}^\perp$ of $a$ and $b$. Then $H = c^\perp$ for a point $c\in \{a,b\}^{\perp\perp}$, by hypothesis $(\ast)$. By Proposition \ref{sub-gen 2}, if $N$ is a generator of $\{a,b\}^\perp$ then $M := \langle N, c\rangle$ is a generator of $\cS$. Clearly $M \subseteq H$, since $H$ contains both $N$ and $c$. This proves $(B)$.

Up till now, we did not use the assumption that $\cS$ is embeddable.
However, to show that $(\ast)$ implies (A), we do assume $(\ast)$ and the embeddability for $\cS$.
We can apply Theorem \ref{Hendrik 1} and find $\cS$ to be symplectic.
Let $\ve:\cS\rightarrow\PG(V)$ be the symplectic embedding of $\cS$ with respect to a non-degenerate alternating form on $V$.
Let $a,b$ be two non-collinear points of $\cS$ and assume $N$ is a generator of $\{a,b\}^\perp$ contained in a generator
$M$ of $\cS$. By Proposition \ref{sub-gen 2}, $N$ is a hyperplane of $M$. Let $m\in M\setminus N$. Then the hyperplane $\overline{H} = \langle \ve(m), \ve(\{a,b\}^{\perp})\rangle_V$  meets
$\langle \ve(\{a,b\}^{\perp\perp}\rangle_V$ in a point $\ve(c)$ with $c\in M\cap\{a,b\}^{\perp\perp}$. Clearly $\overline{H} = \langle \ve(c), \ve(\{a,b\}^\perp)\rangle_V =  \ve(c)^{\perp_\ve}$ and $\ve(M)\subset \overline{H}$. Therefore $c\perp M$. However $M$ is a generator. Hence $c\in M$. So, $c\in M\cap \{a,b\}^{\perp\perp}$.  \hfill $\Box$

\medskip

\noindent
\textbf{The `if' part.} Suppose that both properties $(A)$ and $(B)$ hold for a pair $\{a, b\}$ of non-collinear points of $\cS$. Let $H$ be a hyperplane of $\cS$ and suppose that $H \supseteq \{a,b\}^\perp$. By $(B)$, the hyperplane $H$ contains a generator $M$ of $\cS$ such that $N := M\cap\{a, b\}^\perp$ is a generator of $\{a,b\}^\perp$. By $(A)$, the generator $M$ contains a point $c\in \{a,b\}^{\perp\perp}$. Every line of $\cS$ through $c$ meets $\{a,b\}^\perp$ in a point. Hence $H$ contains $c^\perp$, since it contains $M\cup \{a,b\}^\perp \supseteq \{c\}\cup\{a,b\}^\perp$. It follows that $H = c^\perp$, as claimed in $(\ast)$.  \hfill $\Box$

\subsection{A few comments on Lemma \ref{nuovo}}\label{nuovo sec}\label{MT3}              

In our proof of Lemma \ref{nuovo} the hypothesis that $\cS$ is embeddable has been exploited only to prove that $(\ast)$ implies $(A)$. So, regardless of embeddability, property $(\ast)$ implies $(B)$ for any pair of non-collinear points and, conversely, if $(A)$ and $(B)$ hold for every non-collinear pair of points then $(\ast)$ also holds. One might wonder if the embeddability hypothesis is really necessary in order to prove that $(\ast)$ implies $(A)$. In this section we shall show that indeed it is.

Let $\cS$ be non-embeddable of rank 3. Then all hyperplanes of $\cS$ are singular (Cohen and Shult \cite{CS90}). Hence $(\ast)$ and $(B)$ trivially hold. 

\begin{res}\label{MTA}
Property $(A)$ holds in every non-embeddable polar space of rank $3$, for every non-collinear pair of points.  
\end{res}
\begin{proof} If $\cS$ is the line-grassmannian of a 3-dimensional projective geometry, then the conclusion immediately follows from the fact that every line of $\cS$ belongs to just two planes. Let $\cS$ be a Freudenthal-Tits polar space. For finite rank  property $(A)$ is equivalent to regularity, see Proposition \ref{A0} below. So, proving that all non-collinear pairs of points of $\cS$ are regular is all we need. Actually, this has already been proved in \cite{DBVM1}: Proposition 5.9.4 of \cite{DBVM1} states just this.  \end{proof} 

We shall now discuss two examples of rank 2 where $(\ast)$ holds but $(A)$ fails.     

\begin{example}\label{ABC ex1}
\em
It is well known that the dual of the hermitian generalized quadrangle $H_4(2^2)$ admits no ovoids and no proper full subquadrangles (see e.g. Payne and Thas \cite{PT}). Hence all of its hyperplanes are singular; accordingly, $(\ast)$ trivially holds in it. However its hyperbolic lines have size $2$. Indeed they are just the same as the hyperbolic lines of a (non-full but ideal) subquadrangle isomorphic to $Q_5^-(2)$, which have indeed size 2. Therefore $(A)$ fails to hold in the dual of $H_4(2^2)$.  
\end{example}

\begin{example}\label{ABC ex2}
\em
Let $\cS = P(\cS', x)$ for a generalized quadrangle $\cS'$ of order $q > 3$ and $x$ a regular point of $\cS'$ (notation as in Payne and Thas \cite[\S 3.1.4]{PT}). In view of \cite[\S 1.3.4]{PT}, the quadrangle $\cS$ admits no regular pair of non-collinear points. In other words, since in the finite rank case property $(A)$ is equivalent to regularity, $(A)$ fails to hold in $\cS$ for any non-collinear pair. We shall prove that nevertheless $\cS'$ can be chosen in such a way that $(\ast)$ holds in $\cS$.  

The ovoids of $\cS$ are the sets $O\setminus\{x\}$ for $O$ an ovoid of $\cS'$ containing $x$ (Payne and Thas \cite[\S 3.4.3]{PT}). Suppose that $\cS'$ contains no such ovoid (as when $\cS' = W_3(q)$ with $q$ odd, for instance). Then $\cS$ admits no ovoid. Therefore all non-singular hyperplanes of $\cS$ are subquadrangles. However $q > 3$ by assumption; by \cite[\S 2.2.2(iv)]{PT} we obtain that all full subquadrangles of $\cS$ are grids. The ovoids of a such a grid have size $q$ while the trace of two non-collinear points of $\cS$ has size $q+2$. Hence no non-singular hyperplane of $\cS$ contains the trace of two non-collinear points of $\cS$, namely property $(\ast)$ holds in $\cS$.      
\end{example}

\begin{remark}
When $\cS$ is embeddable its automorphism group acts transitively on the set of non-collinear pairs of points. So, if properties $(A)$ or $(B)$ hold for at least one such pair then they hold for all of them. On the other hand, many non-embeddable generalized quadrangles are known where $(A)$ holds for some non-collinear pairs but not for all of them, namely some but not all pairs of non-collinear points are regular; for instance, this is the case when the quadrangle admits both regular and non-regular points. We believe the same happens for property $(B)$. Regretfully, we have no example at hand to show that this is indeed the case. 
\end{remark}   

\section{Proof of Theorem \ref{Hendrik 1} and Proposition \ref{nuovo H}}\label{sect6}  

\textbf{Proof of Theorem \ref{Hendrik 1}.} The `only if' part of Theorem \ref{Hendrik 1} is well known: let $\cS$ be symplectic and let $\ve:\cS\rightarrow\PG(V)$ be the symplectic embedding of $\cS$. Then $\PG(V) \cong \mathfrak{L}(\cS)$, with $\mathfrak{L}(\cS) = (P, \cL\cup\cL_h)$ defined as at the end of Section \ref{prel}. All singular hyperplanes of $\cS$ lift through $\ve$ to hyperplanes of $\PG(V)$. Hence every singular hyperplane of $\cS$ meets every line of $\mathfrak{L}(\cS)$ non-trivially. Turning to the `if' part, suppose that
\begin{itemize}
\item[$(\dagger)$] every singular hyperplane of $\cS$ meets every hyperbolic line of $\cS$ non-trivially.
\end{itemize} 
\begin{lemma}\label{H1}
For any two points $a, b$ of $\cS$ and a point $c\not\in \{a,b\}^\perp$ there exists a point $x$ such that $x^\perp \supseteq \{a,b\}^\perp\cup\{c\}$.
\end{lemma}
\noindent
\begin{proof} If $a\perp b$ then the unique point of $\langle a, b\rangle\cap c^\perp$ does the job. Suppose that $a\not\perp b$. By $(\dagger)$, $c^\perp$ meets $\{a,b\}^{\perp\perp}$ in a point, say $x$. Clearly, $x^\perp \supset \{a,b\}^\perp\cup\{c\}$.  \end{proof} 

\begin{lemma}\label{H2}
The linear space $\mathfrak{L}(\cS)$ is projective.
\end{lemma}
\noindent
\begin{proof} Let $\cal H$ be the collection of singular hyperplanes of $\cS$. Then $\cal H$ is a collection of geometric hyperplanes of $\mathfrak{L}(\cS)$ (by  $(\dagger)$) and satisfies both the following:
\begin{itemize}
\item[(1)] $\cap_{H\in{\cal H}}H = \emptyset$;
\item[(2)] For any two members $H$ and $H'$ of $\cal H$ a member $H''\in {\cal H}$ exists such that $H'' \supseteq H\cap H'$. 
\end{itemize}
Property (1) trivially follows from the fact that no point of $\cS$ is collinear with all points of $\cS$ while (2) is a way to formulate the conclusion of Lemma \ref{H1}. According to a celebrated theorem of Teirlinck \cite{Teir}, the existence of a family $\cal H$ of hyperplanes satisfying both (1) and (2) characterizes the linear spaces that are projective spaces. Hence $\mathfrak{L}(\cS)$ is projective.  \end{proof} 

Lemma \ref{H2} and Proposition \ref{main thm 1} yield the conclusion: $\cS$ is symplectic.  \hfill $\Box$

\begin{remark}
With the terminology of Hall \cite{H1}, property $(\dagger)$ can also be phrased as follows: the hyperbolic lines of $\cS$ form a {\em copolar space}. In the finite rank case, the statement of Theorem \ref{Hendrik 1} can also be obtained from the classification of copolar spaces by Hall \cite{H1, H2}.
\end{remark}

\noindent
\textbf{Proof of Proposition \ref{nuovo H}.}
The `if' claim is obvious. We prove the `only if' part. Let $H$ be a non-singular hyperplane of a symplectic polar space $\cS$ of rank at least $2$ containing a hyperbolic line $h$ of $\cS$.
Let $\ve$ be the symplectic embedding into the space $V$ equipped with the non-degenerate symplectic form $f$.
Let $h'$ be a hyperbolic line of $\cS$ meeting $h$ in a point $x$ and having a second point $y$ in $H$.
Let $a$ be  the point that maps to the radical of $f$ restricted to the $3$-space $\langle \ve(h),\ve(h')\rangle_V$. As the lines $\langle \ve(a), \ve(y)\rangle$ and $\ve(h)$ are coplanar, $\langle \ve(a), \ve(y)\rangle_V$ meets $\ve(h)$ in a point $\ve(z)$. Clearly, $a \perp z$. Hence $z\perp y$ and $a\in \langle y,z\rangle \subseteq H$. Therefore $a\in H$.
But then all points of $\langle a,h\rangle$ are inside $H$ and, since $h' \subseteq \langle a, h\rangle$, we find $h'$ to be in $H$.

Now take $h'$ to be a hyperbolic line meeting $H$ in at least two points, but not meeting $h$.
If $h'\not\perp h$ then $h'\cap H$ contains a point $x'$ which is not collinear with some point $x$ of $h$. In this case by applying the above first to the hyperbolic lines
$h$ and $h_{x,x'} := \{x,x'\}^{\perp\perp}$ and next $h'$ and $h_{x,x'}$ yields that $h'$ is contained in $H$.
If $h'\perp h$, then take $z$ to be a point off $h$ and $h'$ but on a singular line meeting both $h$ and $h'\cap H$ in points $x$ and $x'$, respectively.
Then $z\in H$.
Fix points $y\in h$ and $y'\in h'\cap H$ different from $x$ and $x'$, respectively.
Applying the above first to the hyperbolic lines $h$ and $h_{y,z}$, next $h_{y,z}$ and $h_{z,y'}$, finally $h_{y',z}$ and $h'$, 
yields that all these hyperbolic lines are in $H$.
This proves $H$ to be a subspace of $\mathfrak{L}(\cS)$.

However $\mathfrak{L}(\cS)$ is a projective space, $\cS$ is a subgeometry of $\mathfrak{L}(\cS)$ and the identity mapping of $\cS$ yields an embedding of $\cS$ in $\mathfrak{L}(\cS)$, isomorphic to the symplectic embedding $\ve$ introduced above. Consequently $\ve(H)$ is a subspace of $\PG(V)$, namely $H$ arises from $\ve$.  \hfill $\Box$

\begin{remark}\label{def trace}
Let $\cS$ be a symplectic generalized quadrangle of infinite order $\mathfrak{s}$. According to Cameron's existence result~\cite{Cam} for ovoids of infinite generalized quadrangles, if $h$ a hyperbolic line of $\cS$ and $X$ a subset of $h$ of cardinality $|X| < \mathfrak{s}$, then $\cS$ admits infinitely many ovoids which contain $X$. However, according to Remark \ref{H rem}, none of them contains $h$. By a slight modification of Cameron's argument, one can also prove that if $|h\setminus X| = \mathfrak{s}$ then ovoids containing $X$ still exist, even if $|X| = \mathfrak{s}$. However none of them contains the whole of $h$. 
\end{remark}

\section{Regular pairs and centric triads}\label{sect5} 

\subsection{Property (A) and Regularity}

Recall that, with the terminology adopted in Sections \ref{Main sec} and \ref{sub-gen}, the maximal singular subspaces of a polar space $\cS$ are its {\em generators} while the singular subspaces which are hyperplanes in some (hence all) of the generators which contain them are the {\em sub-generators} of $\cS$.   

With this terminology, the definition of regular pairs given in Section \ref{Main sec} can be phrased as follows: two non-collinear points $a,b$ of $\cS$ form a \emph{regular  pair} if $\{a,b\}^{\perp\perp}$
equals $\{N,M\}^\perp$ for any pair of disjoint sub-generators $N,M$ of $\cS$ contained in $\{a,b\}^\perp$. If all pairs of non-collinear points of $\cS$ are regular then we say that $\cS$ is {\em regular}. 

The notion of regularity seems to make sense only for finite rank polar spaces.
Indeed, for infinite rank polar spaces the existence of any pair of disjoint  generators  is unclear, see \cite{Pas}. 
However, the following result of \cite{PVM}, already been mentioned in Section \ref{Main sec}, shows we may consider condition (A) to be the right 
generalization of regularity to infinite rank polar spaces:

\begin{prop}\label{A0} Suppose $\cS$ has finite rank. 
A pair of non-collinear points satisfies property $(A)$ if and and only if it is regular.
\end{prop}
\begin{proof}
Suppose $a,b$ is a regular pair. Let $\hat{N}$ be a generator with $N=\{a,b\}^\perp\cap \hat{N}$ a sub-generator.
Then pick a second sub-generator $M$   contained in $\{a,b\}^\perp$,
which is disjoint from $N$. Then for $n\in \hat{N}$ with $n^\perp \supseteq M$, we find $n\in \{a,b\}^{\perp\perp}$.
Notice that $M$ exists in the finite rank case.

This shows that (A) follows from the pair $\{a,b\}$ being regular. Now assume (A) holds for a pair $\{a,b\}$.
Then let $N,M$ be disjoint sub-generators contained in $\{a,b\}^\perp$.
Suppose $c$ is a point in $\{N,M\}^\perp$.
Then $\hat{N}=\langle c,N\rangle$ is a generator of $\cS$, with $N$ as hyperplane. 
So, $\hat{N}$ contains a point $c'$ of $\{a,b\}^{\perp\perp}$.
If $c\neq c'$, then the line $\langle c, c'\rangle$ meets $N$ in a point, which is then in $M^\perp$.
But in a finite rank polar space, we find that $M^\perp\cap N$ is empty. 
We conclude that $c=c'\in \{a,b\}^{\perp\perp}$.
So, $\{N,M\}^\perp=\{a,b\}^{\perp\perp}$ and $\{a,b\}$ is regular. \hfill \end{proof}

The next theorem has been proved in \cite{PIIG}:

\begin{theorem}\label{A1}
Let $\cS$ be embeddable { and of finite rank at least $2$}. Then $\cS$ is regular if and only if $\cS$ admits an embedding $\ve$ of dimension $\dim(\ve) = 2\cdot\rk(\cS)$.
\end{theorem} 

Accordingly, besides symplectic polar spaces, all hermitian varieties of rank $n < \infty$ which admit a $2n$-dimensional embedding are regular and all hyperbolic quadrics are regular (but there is no need of Theorem \ref{A1} to see the latter). 

We shall now turn to a different property. Keeping the assumption that $\rk(\cS) < \infty$, let $M_1, M_2, M_3$ be three mutually disjoint generators of $\cS$. (Such triples of generators always exist except when $\cS$ is either a hyperbolic quadric of odd rank or the line-grassmannian of $\PG(3,\KK)$ with $\KK$ non-commutative.) For $1\leq i < j \leq 3$ let $\pi_{j,i}$ be the orthogonal projection of $M_i$ onto $M_j$ and put $\pi_{1,2.3} := \pi_{1,3}\cdot\pi_{3,2}\cdot\pi_{2,1}$. The bijective mapping $\pi_{1,2,3}$ is a duality of the projective space $M_1$, but in general it is not a polarity; namely in general $\pi_{1,2,3}^2 \neq \mathrm{id}_{M_1}$. The following has been proved in \cite[Theorem 5.12]{PVM}:

\begin{theorem}\label{A2}
The duality $\pi_{1,2,3}$ is a polarity for any choice of three mutually disjoint generators $M_1, M_2$ and $M_3$, if and only if $\cS$ is regular.
\end{theorem} 

We have noticed (Result \ref{MTA}) that Freudenthal-Tits polar spaces are regular. Hence in a Freudenthal-Tits polar space the duality $\pi_{1,2,3}$ is always a polarity. This property is closely related to properties used by Freudenthal  \cite{Fr} to characterize these spaces. See \cite{PVM} for this connection.

\subsection{Property $(B)$ and centric triads}\label{B centric}

Three mutually non-collinear points $a, b$ and $c$ form a {\em centric triad} if $\{a,b,c\}^\perp$ contains a sub-generator of $\cS$. If all triads of mutually non-collinear points of $\cS$ are centric we say that $\cS$ is {\em centric}. 

\begin{lemma}\label{B1}
Let $a, b$ be two non-collinear points of $\cS$ and suppose that the triad $\{a,b,c\}$ is centric for every point $c\not\in \{a,b\}^\perp$. Then $\{a,b\}$ satisfies $(B)$.
\end{lemma}
\begin{proof} Let $H$ be a hyperplane of $\cS$ containing $\{a,b\}^\perp$. The set $\{a,b\}^\perp$ is a subspace but not a hyperplane of $\cS$. Hence $H$ contains a point $c\not\in\{a,b\}^\perp$. By assumption, $\{a, b, c\}^\perp$ contains a sub-generator $N'$ of some generator $N$ of $\cS$. The singular subspace $M := \langle N', c\rangle$  meets $N$ in at least $N'$.
If $M$ is not a generator of $\cS$, it is properly contained in a generator in which we can find a line $\ell$ missing $N'$.
But then take a point $q\in N\setminus N'$. This point is collinear to a point $r\in \ell$. So $\langle N', q\rangle= N \subset r^\perp$,
and hence $N$ is properly contained in the singular space $\langle N,r \rangle$, contradicting
$N$ to be a generator.
Hence also $M$ is a generator of $\cS$. Clearly $M \subseteq H$ and $H$ contains a generator of $\cS$, as required in $(B)$.   \hfill \end{proof}

\begin{theorem}\label{B2}
Let $\cS$ be embeddable. Then $\cS$ is centric if and only if property $(B)$ holds for  every pair of non-collinear points.   
\end{theorem}
\noindent
\begin{proof} The `only if' part is Lemma \ref{B1}. Turning to the `if' part, consider an embedding $\ve:\cS\rightarrow\PG(V)$ of $\cS$. Let $a, b$ and $c$ be three mutually non-collinear points of $\cS$. We know that $\langle \ve(\{a,b\}^\perp)\rangle_V = \{\ve(a),\ve(b)\}^{{\perp_\ve}}$ has codimension $2$ in $\PG(V)$. Clearly, $\ve(c)\not\in \{\ve(a),\ve(b)\}^{{\perp_\ve}}$, since $c\not\in\{a,b\}^\perp$ by assumption and $\ve^{-1}(\{\ve(a),\ve(b)\}^{{\perp_\ve}}) = \{a,b\}^\perp$. Hence $\overline{H} := \langle \ve(c), \{\ve(a),\ve(b)\}^{{\perp_\ve}}\rangle_V$ is a projective hyperplane of $\PG(V)$. Accordingly $H := \ve^{-1}(\overline{H})$ is a geometric hyperplane of $\cS$ and contains $c$ as well as $\{a,b\}^\perp$. In view of $(B)$ on $\{a,b\}^\perp$, the hyperplane $H$ contains a generator $M$ of $\cS$. 

Moreover, since $\{\ve(a), \ve(b)\}^{{\perp_\ve}}$ has codimension 1 in $\overline{H}$, the subspace $\{a,b\}^\perp$ is a geometric hyperplane of $H$. Consequently,  $M$ is also a generator of $H$ and $M\cap\{a,b\}^\perp$ is a sub-generator of $H$, as well as a sub-generator of $\cS$. We can choose the generator $M$ of $H$ in such a way that $c \in M$. With $M$ chosen in this way, $N := M\cap\{a,b\}^\perp$ is a sub-generator of $\cS$ contained in $\{a,b,c\}^\perp$.  \end{proof}

Note that in the previous proof we have exploited only the following seemingly weaker version of $(B)$:  

\begin{itemize}
\item[$(B')$] Let $H$ be a hyperplane of $\cS$ and suppose that $H$ arises from an embedding of $\cS$. If $H\supseteq\{a,b\}^\perp$ then 
$H$ contains a generator of $\cS$.
\end{itemize}
When $\rk(\cS) > 2$ every hyperplane of $\cS$ arises from an embedding of $\cS$ (Shult \cite[7.8.4]{S}).  In this case $(B')$ is just the same as $(B)$. Let $\rk(\cS) = 2$. In this case Theorem 1 of \cite{CGP} implies that all hyperplanes of $\cS$ other than ovoids arise from an embedding of $\cS$. However $\cS$ might admit ovoids which do not arise from any embedding of $\cS$. One might believe that in this case $(B')$ is weaker than $(B)$, but it isn't, as the proof of Theorem \ref{B2} shows. So, the following holds.    

\begin{corollary}\label{B12bis}
Let $\rk(\cS) = 2$ and assume that $\cS$ is embeddable. Suppose that no classical ovoid of $\cS$ contains the trace of two non-collinear points of $\cS$. Then the same holds true for all ovoids of $\cS$, even non-classical ones. In particular, if $\cS$ admits no classical ovoids, then no ovoid of $\cS$ contains the trace of two non-collinear points.
\end{corollary}
The last claim of this corollary follows from the fact that, if $\cS$ admits no classical ovoids, then property $(B')$ is vacuous, hence trivially satisfied. 

Turning back to Theorem \ref{B2}, the hypothesis that $\cS$ is embeddable cannot be dropped from it, as the following examples show. 

\begin{example}
\em 
Let $\cS$ be the dual of $H_4(2^2)$. As noticed in Example \ref{ABC ex1}, property $(B)$ trivially holds in $\cS$. However not all triads of $\cS$ are centric. Indeed let $\ell_1, \ell_2$ be two non-collinear points of $\cS$. Then $\ell_1$ and $\ell_2$ are disjoint lines of ${\cal H} = H_4(2^2)$ and span a hyperplane ${\cal H}'\cong H_3(2^2)$ of $\cal H$. The hyperplane ${\cal H}'$ contains points which do not belong to any of the five lines of ${\cal H}'$ transversal to both $\ell_1$ and $\ell_2$. If $p$ is such a point and $\ell_3$ is a line of $\cal H$ through $p$ not contained in ${\cal H}'$, then $\ell_1, \ell_2$ and $\ell_3$ form a non-centric triad of points of $\cS$.
\end{example} 

\begin{example}
\em
Let $\cS = P(W_3(q), x)$ with $q$ odd, $q > 3$. We know from Example~\ref{ABC ex2} that $(B)$ holds in $\cS$. All triads of mutually non-collinear points of $\cS' := W_3(q)$ are centric and those which are not subsets of hyperbolic lines admit just one center. Let $y_1, y_2, y_3$ be mutually non-collinear points of $\cS$ not in the same hyperbolic line of $\cS'$, let $z$ be the center of $\{y_1, y_2, y_3\}$ in $\cS'$ and suppose that $z\perp x$ in $\cS'$. We shall prove that the triad $\{y_1, y_2, y_3\}$ admits no center in $\cS$.

By way of contradiction suppose that $\{y_1, y_2, y_2\}$ admits a center $z'$ in $\cS$. Clearly $z'\neq z$, as $z\not\in \cS$. It follows that $z'$ is collinear in $\cS'$ with just two of the points $y_1, y_2$ and $y_3$, say $z' \perp y_1, y_2$, and belongs to a hyperbolic line $h$ of $\cS'$ through $x$ and the remaining point $y_3$ of the triad. So, $z \perp x, y_3$ in $\cS'$ and $x, y_3 \in h$. Hence $z\perp z'$ in $\cS'$. However $z, z' \in \{y_1,y_2\}^\perp$ in $\cS'$. We have reached a contradiction. Therefore $\{y_1, y_2, y_3\}$ has no center in $\cS$.  
\end{example}

\begin{example}\label{Bb1}
\em 
Let $\cS$ be the line-grassmannian of $\PG(3,\KK)$ with $\KK$ non-commutative. As noticed in Section \ref{MT3}, property $(B)$ trivially holds in $\cS$. However, no triad of mutually non-collinear points of $\cS$ is centric. Indeed, by way of contradiction, let $a, b, c$ be mutually non-collinear points and suppose that $\{a,b,c\}^\perp$ contains a line. Then $\ell$ belongs to to three distinct singular planes, namely $\langle \ell, a\rangle, \langle \ell, b\rangle$ and $\langle \ell, c\rangle$; this is impossible in this polar space.
\end{example}

\begin{example}\label{Bb2}
\em
Let $\cS$ be a Freudenthal-Tits polar space. We know from Section \ref{MT3} that $(B)$ trivially holds in it. Given two non-collinear points $a, b \in \cS$ and $p\in \{a,b\}^\perp$, let $\cS_p = p^\perp/p$ be the residue of $p$ in $\cS$ and let $\cS^*_p$ be its dual. Then $\alpha := \langle p,a\rangle$ and $\beta := \langle p,b\rangle$ are disjoint lines of $\cS^*_p$. We know that $\cS^*_p$ lives as a quadric in $\PG(11, \KK)$, where $\KK$ is the center of the underlying Cayley algebra of $\cS$ (see Tits~\cite[Chapter 9]{T}, also M\"{u}hlherr~\cite{M} or De Bruyn and Van Maldeghem~\cite{DBVM2}). Accordingly, $\alpha$ and $\beta$ span a grid ${\cal G} \cong Q^+_3(\KK)$ in $\cS^*_p$; the lines of $\cal G$ in the same family as $\alpha$ and $\beta$ are just the lines $\langle p, x\rangle$ of $\cS_p$ for $x\in \{a,b\}^{\perp\perp}$ while the lines $\xi := \langle p, x\rangle$ for $x\in \{a, b, p\}^\perp\setminus\{p\}$ form the other family of lines of $\cG$. However $\cS^*_p$ contains lines disjoint from $\langle \alpha, \beta\rangle$. Let $\gamma$ be one of them, say $\gamma = \langle p, c\rangle$ for a suitable point $c\in p^\perp\setminus\{p\}$. We have $c^\perp\cap (\{a,b\}^{\perp\perp}\cup\{a,b,\}^\perp) = \emptyset$, since $\gamma\cap{\cal G} = \emptyset$. It follows that $\{a,b,c\}^\perp$ contains no lines of $\cS$, namely $a$, $b$ and $c$ form a non-centric triad. 
\end{example}

\vskip.2cm\noindent
\begin{minipage}[t]{\textwidth}
Authors' addresses:
\vskip.4cm\noindent\nobreak
\begin{minipage}[t]{20cm}
\small{Ilaria Cardinali, Antonio Pasini\\
Dep. Information Engineering and Mathematics \\University of Siena\\
Via Roma 56, I-53100 Siena, Italy\\
ilaria.cardinali@unisi.it,  antonio.pasini@unisi.it}
\end{minipage}
\vskip.5cm\noindent\nobreak
\begin{minipage}[t]{20cm}
\small{Hans Cuypers,\\
Dep. Mathematics and Computer Sciences,\\
Eindhoven University of Technology,\\
5600 MB Eindhoven, The Netherlands}\\
f.g.m.t.cuypers@tue.nl
\end{minipage}
\vskip.5cm\noindent\nobreak
\begin{minipage}[t]{20cm}
\small{Luca Giuzzi\\
D.I.C.A.T.A.M. \\
University of Brescia\\
Via Branze 43, I-25123 Brescia, Italy \\
luca.giuzzi@unibs.it}
\end{minipage}
\end{minipage}


\begin{thebibliography}{999}

\bibitem{BC} F. Buekenhout and A. M. Cohen. {\em Diagram Geometry,} Springer, Heidelberg 2013.

\bibitem{Cam} P. J. Cameron,  Ovoids in infinite incidence structures, {\em Arch. Math.} {\bf  62} (1994), 189-192.

\bibitem{CGP} I. Cardinali, L. Giuzzi and A. Pasini. Nearly all subspaces of a classical polar space arise from its universal embedding, {\em Lin. Alg. Appl.} {\bf  627} ( 2021), 287-307.

\bibitem{CS90} A. M. Cohen and E. E. Shult. Affine polar spaces, {\em Geom. Dedicata} {\bfseries 35} (1990), 43-76.


\bibitem{CJP} H. Cuypers, P. Johnson and A. Pasini. On the embeddability of polar spaces, {\em Geom. Dedicata} {\bf 44} (1992), 349-358.

\bibitem{DB} B. De Bruyn and A. Pasini. On Symplectic Polar Spaces over non-Perfect Fields of Characteristic $2$, {\em Lin. Multilinear Alg.} {\bf 57} (2009), 567-575.

\bibitem{DBVM1} B. De Bruyn and H. Van Maldeghem.  Non embeddable polar spaces, {\em M\"{u}nster J. Math.} {\bf 7} (2014), 557-588.

\bibitem{DBVM2}  B. De Bruyn and H. Van Maldeghem. Universal and homogeneous embeddings of dual polar spaces of rank 3 defined over a quadratic alternative division ring, {\em J. Reine Angew. Math.} {\bf 715} (2016), 39-74.

\bibitem{Fr} H. Freudenthal. Beziehung der $E_7$ und $E_8$ zur Oktavenebene, I-XI, {\em Proc. Kon. Ned. Akad. Wet. Ser. A}, {\bf 57} (1954), 218-230, 363-368: {\bf 57} (1955), 151-157, 277-285; {\bf 62} (1959), 165-201, 447-474; {\bf 66} (1963), 457-487.

\bibitem{H1} J. I. Hall. Classifying copolar spaces and graphs, {\em Quart. J. Math.}, {\bf 33} (1982), 421-449. 

\bibitem{H2} J. I. Hall. Identifying classical geometries, in {\em Finite Geometries}, vol. 82 of {\em Lecture Notes Pure and Appl. Math.}, Dekker, New York (1983), 175-195. 

\bibitem{J1} P. M. Johnson. Polar spaces of arbitrary rank, {\em Geom. Dedicata} {\bf 35} (1990), 229-250. 

\bibitem{J2} P. M. Johnson. Semiquadratic sets and embedded polar spacs, {\em J. Geometry} {\bf 64} (1999), 102-127. 

\bibitem{M} B. M\"{u}hlherr. A geometric approach to non-embeddable polar spaces of rank 3, {\em Bull. Soc. Math. Belgique} {\bf 42} (1990), 577-594.


\bibitem{Pas} A. Pasini. On polar spaces of infinite rank, {\em J. Geometry} {\bf 91} (2008), 84-118.


\bibitem{PIIG} A. Pasini. Synthetic and projective properties of embeddable polar spaces, submitted to {\em Inn. Inc. Geometry}. 

\bibitem{PVM} A. Pasini and H. Van Maldeghem. An essay on Freudenthal-Tits polar spaces, submitted to {\em J. Algebra}.  

\bibitem{PT} S. E. Payne and J. A. Thas. {\em Finite generalized quadrangles.} Second Edition. EMS Series of Lectures in Mathematics, European Mathematical Society (EMS), Z\"{u}rich, 2009.




\bibitem{S} E. E. Shult. {\em Points and Lines,} Springer, Berlin, 2010.

\bibitem{Teir} L. Teirlinck. On projective and affine hyperplanes, {\em J. Combin. Th. Ser. A}, {\bf 28} (1980), 290-306. 

\bibitem{T} J. Tits. {\em Buildings of Spherical Type and Finite BN-Pairs}, Springer Lecture Notes in Math. {\bfseries 386} (1974).


\bibitem{HVM1} H. Van Maldeghem. {\em Generalized Polygons}, Birkh\"{a}user, Basel 1998.

\bibitem{HVM2} H. Van Maldeghem, private communication to the fourth author.

\end{thebibliography}
\end{document}